\documentclass[12pt]{article}

\usepackage{amsmath}
\usepackage{amsfonts}
\usepackage{amssymb}
\usepackage{graphicx, color}
\usepackage{hyperref}
\usepackage{epsfig}
\usepackage{epstopdf} 
\definecolor{darkgreen}{rgb}{0,0.55,0}

\newtheorem{proposition}{Proposition}[section]
\newtheorem{theorem}{Theorem}[section]

\newtheorem{remark}[theorem]{Remark}

\DeclareSymbolFont{AMSb}{U}{msb}{m}{n}
\DeclareMathSymbol{\N}{\mathbin}{AMSb}{"4E}
\DeclareMathSymbol{\Z}{\mathbin}{AMSb}{"5A}
\DeclareMathSymbol{\R}{\mathbin}{AMSb}{"52}
\DeclareMathSymbol{\Q}{\mathbin}{AMSb}{"51}
\DeclareMathSymbol{\I}{\mathbin}{AMSb}{"49}

\newcommand{\Si}{{\mathbb S}}
\newcommand{\Ce}{{\mathbb C}}

\newcommand{\calK}{{\mathcal K}}
\newcommand{\calM}{{\mathcal M}}

\newcommand{\calH}{{\mathcal H}}

\newcommand{\calS}{{\mathcal S}}
\newcommand{\calT}{{\mathcal T}}

\numberwithin{equation}{section}

\begin{document}

\title{A Sharp Inequality on the Exponentiation of Functions  on the Sphere}

\author{{Sun-Yung Alice Chang }
\qquad Changfeng Gui
}
\date{\today}

\smallbreak \maketitle

\begin{abstract}
In this paper we show  a new inequality which generalizes to the unit  sphere the Lebedev-Milin inequality  of the exponentiation  of functions on the unit circle.  It may also be regarded as  the counterpart   on the sphere   of the second inequality in  the Szeg\"o limit theorem on the Toeplitz determinants on the circle.  On the other hand,  this inequality  is also   a variant of several classical inequalities of Moser-Trudinger type on the sphere. The inequality incorporates  the deviation of the  center of mass from  the origin  into the optimal inequality of Aubin  for functions with  mass centered  at the origin,  and improves  Onofri's  inequality with the contribution of the shifting of the  mass center  explicitly expressed.

\end{abstract}
\maketitle

\section{Introduction} 

Let $\Si^2$ be the unit sphere and for $u\in H^1(\Si^2)$ define 

\begin{equation}\label{MTAO}
F_\alpha(u)=\alpha \int_{\Si^2}|\nabla u| ^2 d\omega+ 2\int_{\Si^2}u d\omega-\log \int_{\Si^2}e^{2u}d \omega,
\end{equation}
where the volume form $d\omega$ is normalized so that $\int_{\Si^2}d\omega=1$. The well-known Moser-Trudinger inequality \cite{Moser-MR0301504} says that $F_\alpha$ is bounded below if and only if $\alpha \geq 1$. 
Later on Onofri \cite{O-MR677001} sharpened Moser-Trudinger inequality and showed that for $\alpha \geq 1$ the best lower bound of $F_\alpha$ is equal to zero. Onofri's inequality was based on an inequality established earlier by  Aubin \cite{Aubin2-MR534672} who proved that if $F_\alpha$ is restricted to  
\[\mathcal{M}:=\{u\in H^1(\Si^2): \ \ \int_{\Si^2}e^{2u}x_i=0, \ \ i=1,2,3\},\]
then for $\alpha >\frac{1}{2}$, $F_\alpha$ is bounded below and the infimum is attained in $\mathcal{M}$. All these inequalities play crucial roles in the ``Nirenberg's problem" of prescribing Gaussian curvature, 
in particular in the work of Chang and Yang 
(\cite{CY1-MR908146} and  \cite{CY2-MR925123}).  In their effort to prescribe Gaussian curvature without additional assumption on the symmetry of the curvature,  Chang and Yang have further improved the above Aubin-Onfri inequality by showing that (see Proposition B  in \cite{CY1-MR908146})  for $\alpha$ sufficient close but less than 1, the lower bound of $F_\alpha $ again is equal to zero for $u$ in the class $\mathcal{M}$,  their work led to the following conjecture: \\ \\
{\bf Conjecture A.} For $\alpha \geq \frac{1}{2}$ 
\[\inf_{u\in \mathcal{M}} F_{\alpha}(u)=0.\]

In 1998,  Feldman, Froese, Ghoussoub and Gui \cite{FFGG-MR1606461} proved that this conjecture is true for axially symmetric functions when $\alpha>\frac{16}{25}-\epsilon$. Later the second author and Wei \cite{GW-MR1760786}, and independently Lin \cite{Lin1-MR1770683} proved Conjecture A for axially symmetric functions. In \cite{GL-MR2670931}  Ghoussoub and Lin showed that Conjecture A holds true for $\alpha \geq \frac{2}{3}-\epsilon$, for some $\epsilon>0$. Finally  Gui and Moradifam  proved in \cite{GM1}  that  Conjecture A is indeed true.  Actually they \cite{GM1} obtained something stronger than the conjecture, by showing  the following uniqueness result for  the corresponding Euler-Lagrange equation for the functional $F_\alpha$. 

\begin{theorem}\label{MAOTheorem}
 The following  equation 
\begin{equation}\label{standardPDE}
\alpha\Delta u+\frac{e^{2u}}{\int_{\Si^2}e^{2u}  d\omega}-1=0 \ \ \hbox{on}\ \ \Si^2
\end{equation}
has only constant solutions for $\frac{1}{2} \le \alpha <1$. 
\end{theorem}

\section{A Refined Aubin-Onofri Type Inequality} 

The main result in this paper is to establish a variant of Aubin-Onofri inequality. To motivate the study of such type of inequalities, we first recall the  classical Lebedev-Milin inequality on the exponentiation  of functions  defined on the unit circle $\Si^1$, which is in spirit similar to that of the 
Moser-Trudinger inequality for functions defined on $\Si^2$. 

Assume on $\Si^1 \subset \R^2 \sim  \Ce$
$$
u(z)= \sum_{k=1}^\infty a_k z^k, \quad e^{u(z)}= \sum_{k=0}^\infty \beta_k z^k.
$$
Then  the Lebedev-Milin inequality on  the unit circle  (\cite{D}) states 
\begin{equation}\label{LM}
\log(\sum_{k=0}^\infty |\beta_k|^2 ) \le  \sum_{k=1}^\infty k |a_k|^2
\end{equation}
if the right hand side is finite , and equality holds if and only if $ a_k=\gamma^k/k$ for some $\gamma \in \Ce$ with $|\gamma|<1$.   This is well known in the  community of  univalent functions,  in particular in connection with Bieberbach conjecture.

Denote $D$ the unit disc on $\mathbb R^2$. For any real  function $u$ defined on the unit circle, we recall  that the right hand side of \eqref{LM} is indeed  $H^{\frac{1}{2} }(\Si^1)$  norm of $u$,   which can also be  identified as the  $H^1 (D) $ norm of the
harmonic extension, which we denote again by $u$,  on the disc $D$.  Then the classical Lebedev-Milin inequality  may be written as  
\begin{equation}\label{ML}
\log (\frac{1}{2\pi} {\int_{\Si^1}  e^u d \theta} ) -  \frac{1}{2 \pi} { \int_{\Si^1} u d\theta } \leq \frac{1}{4\pi} || \nabla u||^2_{L^2(D)} .
\end{equation}
It turns out Lebedev-Milin inequality is the ``first step'' of a string of monotonically increasing inequalities in the  Szeg\"o Limit Theorem (\cite{GS-1958}, 5.5a) on Toeplitz determinants. 
Here we will just quote the second inequality in the Szeg\"o limit theorem:
\begin{equation} \label{GS2}
\log  (|\frac{1}{2\pi} {\int_{\Si^1} e^u d \theta} |^2 -  |\frac{1}{ 2\pi} {\int_{\Si^1} e^u e^{i\theta}d \theta }|^2 ) -  \frac{1}{ \pi} {\int_{\Si^1} u d\theta } \leq \frac {1}{4\pi} || \nabla u||^2_{L^2(D)}.
\end{equation}
One notes that  in the special case when $ {\int_{\Si^1} e^u e^{i\theta}d \theta }=0$, as a direct consequence of the above inequality we have  
\begin{equation} \label{GS2cor} 
\log (\frac{1}{2 \pi} {\int_{\Si^1} e^u d \theta} ) -  \frac{1}{2 \pi} {\int_{\Si^1} u d\theta}  \leq \frac {1}{8\pi} || \nabla u||^2_{L^2(D)}.
\end{equation}
Indeed this special form of the inequality was independently verified by Osgood, Phillips, Sarnak \cite{OPS88} and was used in their study of isospectral compactness for metrics defined on compact surfaces.  It was later pointed 
out by H. Widom (\cite{W88}) that it is a direct consequence of the Szeg\"o Limit Theorem.  We remark that actually Widom has also pointed out that for all integer $k$,  there is a string of such inequalities for functions $u$ with $\int_{\Si^1} e^u e^{ij\theta} d\theta = 0$
for all $ 1 \leq j\leq k$. In a recent work(\cite{CH2019}), Chang and Hang have further explored this angle and established a weaker form of such inequalities for functions defined on the 2-sphere with vanishing higher order of moments. 

The relevance to us is the apparent comparison of the inequality of \eqref{GS2cor} on $\Si^1$ as compared to Conjecture A 
in the introduction for functions defined  on $\Si^2$. This leads us to ask the question if there a corresponding
inequality on $\Si^2$ similar to that of \eqref{GS2}, which in the special case when $u$ is in $\mathcal{M}$  reduces to the statement in Conjecture A. 
  
Motivated by this, we consider the following family of  functionals in $H^1(\Si^2)$: 
\begin{equation}\label{I_alpha}
I_\alpha(u)=\alpha \int_{\Si^2}|\nabla u| ^2 d\omega+ 2\int_{\Si^2}u d\omega-\frac{1}{2} \log [(\int_{\Si^2}e^{2u}d \omega)^2- \sum_{i=1}^{3}(\int_{\Si^2}e^{2u} x_i d \omega)^2 ]
\end{equation}
where $\alpha >0$. 

The question we are asking is what is the minimum value of $\alpha$ for which the functional $F_\alpha (u)$ stays non-negative for all functions $ u \in H^{1}(\Si^2) $.  One notices that if such a minimum value $ \alpha$ is  $\frac{1}{2} $, then we would 
recover the statement in  Conjecture A.  But to our surprise, the answer of the question is actually no, and the minimum value of such $ \alpha $ is actually $ \frac{2}{3} $. We will present here our analysis, and state the following result 
as our main theorem.

\begin{theorem} \label{main}
For  any $\alpha >0$,  we have 
\begin{equation}\label{ineq}
I_{\alpha} (u) \ge (\alpha-\frac{2}{3})  \int_{\Si^2}|\nabla u| ^2 d\omega, \quad \forall u \in H^1(\Si^2).
\end{equation}
In particular, when   $\alpha \ge \frac{2}{3}$  we have 
\begin{equation}\label{ineq1}
I_{\alpha} (u) \ge 0, \quad \forall u \in H^1(\Si^2).
\end{equation}
Furthermore,  for $ 0 <  \alpha < \frac{2}{3} $,  $ \inf_{H^1(\Si^2)} I_{\alpha} (u)  = - \infty$.
\end{theorem}

In the rest of the section, we will present the proof of the above theorem. Due to the invariance of $I_{\alpha} (u)$ by a constant addition, we may confine our discussion in the normalized space
\begin{equation}\label{normalized}
\calH= \{ u \in H^1(\Si^2):  \int_{\Si^2}e^{2u}d \omega=1\}. 
\end{equation}

 The strategy is to first study the Euler-Lagrange equation of the functional $I_{\alpha}$, assuming the critical point is obtained. It turns out for the special value $\alpha = \frac {2}{3}$, for each
point $\vec a$ in the unit ball $B_1\subset \R^3$, there is a unique solution $u\in \calH$, which we can write down explicitly,  of the Euler-Lagrange equation, with the center of  mass of $e^{2u} $ being at $\vec{a} $.   Based on this analysis, we then study the minimum of $I_{\alpha}(u)  $ over the class of $u$ with a fixed center of mass and verify  that it is achieved for each $ \alpha > \frac{1}{2} $.  Although  the infimum of $I_{\alpha}(u) $  tends to negative infinity  as $\vec a$ goes to the unit sphere $\Si^2=\partial B_1$  when $ \frac {1}{2} < \alpha < \frac{2}{3} $,  it turns out that $I_{\alpha} (u) \geq 0$ for all $u \in H^{1} (\Si^2)$  when $ \alpha \geq \frac{2}{3} $,  due to the complete understanding of the
critical points of  $I_{2/3} (u)$  in $\calH$ with  the  center  constrained. 
\vskip .2in 

We now begin the analysis.  For each $ u \in \calH$, denote
\begin{equation}\label{a_i}
a_i=\int_{\Si^2}e^{2u} x_i d \omega, \quad  i=1, 2, 3.
\end{equation}

\begin{proposition} \label{euler} 
The Euler Lagrange equation for the functional $I_{\alpha}$ in $\calH$  is
\begin{equation}\label{simple}
\alpha \Delta u+ \frac{ 1 -\sum_{i=1}^{3} a_i x_i }{1-\sum_{i=1}^{3}a_i^2} e^{2u}-1=0 \ \ \hbox{on}\ \ \Si^2.
\end{equation}
\end{proposition}

%
 


We now study the solution of equation \eqref{simple}.

\begin{proposition}\label{eq-main} 

i )  When $\alpha \in (0, 1) $ and $ \alpha \not = \frac{2}{3}$,  equation \eqref{simple} has only zero solution in $\calH$; 

ii) When $\alpha =\frac{2}{3}$,  for any $\vec{a}=(a_1, a_2, a_3) \in B_1$, there is a unique solution $u$  to equation \eqref{simple} in $\calH$  such that  \eqref{a_i} holds.   In particular,  $u$ is axially symmetric about $\vec{a}$  if $\vec{a} \not =(0,0,0)$.   After a proper rotation,  the solution  $u$ is explicitly given by  the formula in \eqref{solution} below.
 
\end{proposition} 

\noindent {\bf Proof} :
\vskip .1in


To investigate \eqref{simple},  recall  the Kazdan-Warner condition for the Gaussian curvature equation:
\begin{equation}\label{ageneral}
\Delta u +K (x) e^{2u}= 1 \ \ \hbox{on}\ \ \Si^2,
\end{equation}
then 
\begin{equation} \label{KW} 
\int_{\Si^2} ( \nabla K (x) \cdot  \nabla x_j ) e^{2u} d \omega =0 \,\,\,\, \hbox{for each j= 1,2, 3}. 
\end{equation}
If $u$ satisfies the \eqref{simple}, then

\begin{equation}\label{gaussian}
K (x) = \frac{1}{\alpha} \frac{ (1 -\sum_{i=1}^{3} a_i x_i ) }{(1-\sum_{i=1}^{3}a_i^2) }  + ( 1 - \frac{1}{\alpha} ) e^{-2u}.
\end{equation}

Substituting \eqref{gaussian} into \eqref{KW}, we obtain for each $j =1,2,3$, 
$$
\frac{1}{\alpha} \frac {1}{ (1-\sum_{i=1}^{3}a_i^2)} \int_{\Si^2} (\sum_i a_i  (\nabla x_i  \cdot \nabla  x_j))  e^{2u} d \omega \, =\,  (-2 ) ( 1- \frac{1}.{\alpha} ) \int_{\Si^2} (\nabla u  \cdot \nabla x_j) d\omega.
$$
Integrate by part the last term, and substitute the term $\Delta u$ in equation \eqref{simple} and simplify.  We then  get
$$
\sum_i a_i \int_{\Si^2} \nabla x_i  \cdot \nabla x_j e^{2u} d \omega \, = \,  2 ( 1- \frac{1}{\alpha} ) \sum_i a_i \int_{\Si^2} x_i x_j e^{2u}  d\omega -2 ( 1- \frac{1}{\alpha} ) a_j.
$$
We now notice that $\nabla x_i \cdot \nabla x_j = - x_i x_j$ when $ i \not = j$, $ |\nabla x_j |^2 = 1- x_j^2$, thus we get 
$$
 ( 3- \frac{2}{\alpha} ) a_j \int_{\Si^2} e^{2u}  d\omega = ( 3- \frac{2}{\alpha} ) \sum_i a_i \int_{\Si^2} x_i x_j e^{2u} \omega.
$$
Multiply the above formula by $a_j$ and sum over $j$, we get 
$$
 ( 3- \frac{2}{\alpha}) \int_{\Si^2} \bigl( \sum_j a_j^2- |\sum_i a_i x_i|^2 \bigr)  e^{2u} d\omega=0.
$$

This implies that if $\alpha \not = \frac{2}{3}$,  there holds $a_i=0, i=1, 2, 3$,    since
$$
(\sum_{i=1}^{3} a_i x_i )^2  \le \sum_{i=1}^{3} a_i^2  \ \ \hbox{on}\ \ \Si^2
$$
and the equality  only holds  when  $ \vec{a}=(a_1, a_1, a_3)$ is the  zero vector or $ x$ is parallel to  the vector $(a_1, a_1, a_3)$ if it is not the zero  vector. 
Therefore,  we conclude that when $\alpha \in (0, \frac{2}{3}) \cup (\frac{2}{3}, 1)$,  the equation \eqref{simple}  have only zero solution, in view of Theorem \ref{MAOTheorem}.

For $\alpha= \frac{2}{3}$,  we assume that $u$ is a solution to the coupled equations \eqref{a_i} and \eqref{simple}. Without of loss of generality,  we may assume that $(a_1, a_2, a_3)=(0, 0, a)$ 
with $a \in (0, 1)$ and consider 
 \begin{equation}\label{eq-a}
\frac{2}{3}  \Delta u+ \frac{ 1 -a x_3} {1-a^2} e^{2u}-1=0 \ \ \hbox{on}\ \ \Si^2.
\end{equation}
We shall use the stereographic projection to transform the equation to be on  $\R^2$. 
Let $\Pi$ be the stereographic projection $\Si^2\rightarrow \R^2$ with respect to the north pole $N=(0,0, 1)$: 
\[y=\Pi(x):= \left( \frac{x_1}{1-x_3}, \frac{x_2}{1-x_3}\right).\]
Note that 
$$
x_3=\frac{|y|^2-1}{|y|^2+1},   \quad  d\omega=\frac{dy}{\pi (|y|^2+1)^2}  
$$
Suppose $u$ is a solution of \eqref{eq-a}, and let
\[w(y):=u(\Pi^{-1}(y))-\frac{3}{2}  \ln(1+|y|^2 ) \ \ \hbox{for}\ \ y\in \R^2.\]
Then $w$ satisfies 
\begin{equation}\label{plane}
\Delta w+\frac{6}{1+a} (\mu^2+|y|^2)e^{2w}=0 \ \ \hbox{in} \ \ \R^2
\end{equation}
where $\mu^2=\frac{1+a}{1-a}>1,  b>0$ 
and 
\begin{equation}\label{total}
\int_{\R^2}(\mu^2+|y|^2)e^{2w} dy= (1+a) \pi. 
\end{equation}
Now it is easy to verify directly  that
$$
w(y)= -\frac{3}{2} \ln(\mu^2+|y|^2)+2 \ln \mu+\frac{1}{2} \ln  \frac{2}{1+\mu^2}
$$
is a solution to \eqref{plane} and \eqref{total}, and hence   $u(x)$  defined by
\begin{equation}\label{solution}
u(x)=u(\Pi^{-1}(y)): =\frac{3}{2} \ln\frac{1+|y|^2 }{\mu^2+ |y|^2} +2\ln \mu +\frac{1}{2} \ln  \frac{2}{1+\mu^2}
\end{equation}
is a solution to \eqref{eq-a}.
It is also easy to compute that  $ \int_{\Si^2} e^{2u} d\omega=1$ and  
$$
\int_{\Si^2} e^{2u} x_3d\omega=\int_{\R^2}e^{2u(\Pi^{-1}(y))} ( \frac{|y|^2-1}{|y|^2+1}) \frac{dy}{\pi (|y|^2+1)^2} = a,
$$
and therefore $u$ is a solution to \eqref{a_i} and \eqref{simple} with $(a_1, a_2, a_3)=(0, 0, a)$.

To show the uniqueness of the solution to \eqref{eq-a},   we will recall a general result regarding the  radial symmetry of solutions. 
Assume  $u \in C^2 (\R^2) $  satisfies 
\begin{equation}\label{general}
\Delta u+\calK(|y|)e^{2u}=0 \ \ \hbox{in}\ \ \R^2,
\end{equation}
and 
\begin{equation}\label{16Pi}
\frac{1}{2\pi} \int_{\R^2}\calK(|y|)e^{2u}dy =\beta < \infty,
\end{equation} 
where  $\calK(y):=\calK(|y|) \in C^2(\R^2)$ is  a non constant  positive function  satisfying 
\begin{eqnarray*}
&(K1) \quad \quad   &\Delta \ln(\calK(|y|))\geq 0,  \quad y \in \R^2\\
&(K2)  \quad \quad  & \lim_{|y| \to \infty} \frac{|y| \calK'(|y|)}{\calK(|y|)}=2l>0 ,  \quad y \in \R^2.
\end{eqnarray*}

The following general symmetry result is proven in \cite{GM1}.

\begin{proposition} \label{th-general}
Assume that $\calK(y)=\calK(|y|)>0$  satisfies  $(K1)-(K2),$  and  $u$ is a solution to \eqref{general}-\eqref{16Pi}  with $ l+1< \beta \le  4$.  
Then $u$ must be  radially symmetric. 
\end{proposition}

Applying  Proposition \ref{th-general} to \eqref{plane} and  \eqref{total} with $l=1, \beta=3$,  we conclude that the solution to \eqref{plane} and  \eqref{total}  must be radially symmetric.  Furthermore,   such a radial solution must be unique by Theorem 1.5 of \cite{Lin1-MR1770683}. 
Therefore we have finished the proof of the Proposition \eqref{eq-main}.  $\blacksquare$
\vskip .1in
%
%
 %
  

\noindent{\bf Proof of Theorem \ref{main}}
\vskip .1in

For any $\vec{a}=(a_1, a_2, a_3) \in B_1:=\{ |a|<1\} \subset \R^2 $,  let us define
\begin{equation}\label{minimizing}
\calM_{\vec{a}}:=\{u\in \calH \subset H^1(\Si^2): \ \ \int_{\Si^2}e^{2u}x_i=a_i, \ \ i=1,2,3\}.
\end{equation}

First we consider a  constrained minimizing problem  on $\calM_{\vec{a}}:$
$$
m(\alpha, \vec{a}):=\min_{u \in \calM_{\vec{a}} }  I_{\alpha} (u)
$$
and recall the following compactness result:

\begin{proposition} \label{compact}
For  any $\alpha >\frac{1}{2}, \vec{a}=(a_1, a_2, a_3) \in B_1$,  there exists $C_{\alpha, \vec{a}} \in \R $  such that
\begin{equation}\label{lowerbound}
I_{\alpha} (u) \ge C_{\alpha, \vec{a} },  \quad \forall u \in \calM_{\vec{a}}.
\end{equation}
Furthermore,   there is a positive  constant  $M_{\alpha, |\vec{a}|, C} >0$  depending only on $ \alpha, |\vec{a}|<1$ and $C$  such that $ \| u\|_{H^1(\Si^2)}  \le M_{\alpha, |\vec{a}|, C} $  in the sub level set  $ I_{\alpha, \vec{a} }^C: =\{ u \in \calM_{\vec{a} },  I_{\alpha} (u) \le C  \}$.
\end{proposition}

\noindent {\bf Proof}:
\vskip .1in

This result may be known to researchers in the area,  although it seems not stated or proven explicitly  in the literature.  Here we will give a sketch of  proof following  Proposition 2.1 of  \cite{CH2019}. 

Assume that for some $\alpha >\frac{1}{2},  \vec{a}=(a_1, a_2, a_3) \in B_1$,  there is a sequence 
 $u_k \in \calM_{\vec{a}}, k=1, 2, \cdots $ such that 
 $I_{\alpha} (u_k)  \to -\infty$ as $k \to \infty$.
 Then
 $$
 {\bar u}_k : =  \int_{\Si^2} u_k d\omega  \to -\infty, \quad k \to \infty.
 $$
 By the classical Moser-Trudinger inequality, we have
$$
 \int_{\Si^2}|\nabla u_k| ^2 d\omega\ge - 2{\bar u}_k  \to \infty,  \quad k \to \infty.
$$
 Let  $m_k= (\int_{\Si^2}|\nabla u_k| ^2 d\omega)^{\frac{1}{2}}$  and
 $$
   v_k =\frac{u_k- {\bar u}_k} {m_k}. 
 $$ 
 Then,  when $k$ is sufficiently large,  by the assumption $v_k$ satisfies 
 $$
 \ln( \int_{\Si^2}  e^{2m_k v_k} d \omega)  = -2\bar{u}_k \ge \alpha \int_{\Si^2}|\nabla u_k| ^2 d\omega =\alpha m_k^2.
  $$
 Assume that  $v_k$ converges weakly  to $v$ in $H^1(\Si^2)$ as  $ k \to \infty$
 and 
 $$
 |\nabla v_k| ^2 d\omega  \to  |\nabla v| ^2 d\omega +\sigma  \quad  \hbox{and } \quad \frac{ e^{2m_k v_k} d \omega}{ \int_{\Si^2}  e^{2m_k v_k} d \omega} \to \nu, \quad k \to \infty.
$$
 in measure, where $\sigma(\Si^2)=\nu (\Si^2)=1$.  
  Then,  by  Proposition 2.1 in \cite{CH2019}  we have, 
 $$
 \{x \in \Si^2:  \sigma(x)  \ge \alpha \} =\{P\}  \quad \hbox{and } \quad  \nu=\delta_{P}
 $$ 
  for some $P \in \Si^2$ since $\alpha >\frac{1}{2}$.  Note that here we have a normalized area of the unit sphere being 1 with the measure $\omega$ while in  \cite{CH2019}   the area of the unit sphere is $4 \pi$. 
  
 This leads to a contradiction that
 $$
\vec{a}= \int_{\Si^2} e^{2u_k} x d \omega = \frac{  \int_{\Si^2} e^{2m_k v_k} x d \omega}{ \int_{\Si^2}  e^{2m_k v_k} d \omega }  \to P, \quad k \to \infty.
$$
The argument also shows that $\| \nabla u\|_{L^2(\Si^2)} $ is bounded in  the sub level set $I_{\alpha, \vec{a}}^C:=\{ u \in \calM_{\vec{a}},  I_{\alpha} (u) \le C \}$ for any fixed $C \in \R$.  The  Jensen inequality and the convexity of the exponential function  as well as \eqref{lowerbound} lead to the boundedness of  ${\bar u}=\int_{\Si^2} u d \omega $  in  the set $ I_{\alpha, \vec{a} }^C$.

Therefore, Proposition \ref{compact} holds.  $\blacksquare$
\vskip .1in
 
From here it is standard to show that there exists a  minimizer  $u_{\alpha, \vec{a}} \in \calM_{\vec{a}} $ of \eqref{minimizing}  satisfying 
 
\begin{equation}\label{minimizer}
\alpha \Delta u+e^{2u}(\rho-\sum_{i=1}^3 \beta_i x_i)=1, \quad x \in \Si^2
\end{equation}
for some  $\rho \in \R$ and $\vec{\beta}=(\beta_1, \beta_2, \beta_3) \in \R^3$.

To be more precise, for a fixed $\alpha >\frac{1}{2}$ there is a minimizing sequence of $u_k \in \calM_{\vec{a}}, k=1, 2, \cdots $  of  $I_{\alpha} $ such that  $u_k$ is bounded in $H^1(\Si^2)$ 
and  $ u_k  $  converges weakly to $u_{\alpha, \vec{a}}$ in $H^1(\Si^2)$ and  $u_{\alpha, \vec{a}} \in \calM_{\vec{a}}$. (See, e.g., the proof of Theorem 5.1  of \cite{CH2019}.)    Hence $u_{\alpha, \vec{a}}$ is a minimizer of  $\min_{u \in \calM_{\vec{a}} }  I_{\alpha} (u). $  
It is easy to see that
$$
\rho=1 + \sum_{i=1}^3 \beta_i  a_i. 
$$

Using Kazdan-Warner condition \eqref{KW},  we obtain 

\begin{equation}\label{condition}
2 (\frac{1}{\alpha} -\frac{3}{2})  \sum_{i=1}^3 \beta_i \int_{\Si^2} x_ix_j e^{2u} d\omega =
 2( \frac{1}{\alpha}-1)\rho  a_j-\beta_j, \quad j=1, 2, 3.
 \end{equation}
 
 In particular, when $\alpha=\frac{2}{3}$,  we have
 $$
 \beta_j= \frac{a_j}{1-|\vec{a}|^2}, \quad j=1, 2, 3.
 $$
 
 Then equation \eqref{minimizer} is equivalent to \eqref{simple} when $\alpha=\frac{2}{3}$. 
 
 After a proper rotation so that $\vec{a}$ points  to the north pole and using the  stereographic project $\Pi: \Si^2 \to \R^2$,  the solution is uniquely determined by 
 $$
 u_{\frac{2}{3}, |\vec{a}|} (x):=\frac{3}{2} \ln\frac{1+|y|^2 }{\mu^2+ |y|^2} +2\ln \mu +\frac{1}{2} \ln  \frac{2}{1+\mu^2}
$$
where $\mu^2=\frac{1+|\vec{a}|}{1-|\vec{a}|}>1$.

Hence, by direct computations we have
\begin{equation}\label{gradient}
\begin{split}
\int_{\Si^2} |\nabla  u_{\frac{2}{3}, |\vec{a}|} |^2 d\omega  
&= \frac{1}{4\pi}  \int_{\R^2} \frac{9 (\mu^2-1)^2}{4} \frac{|y|^2}{(|y|^2+1)^2(|y|^2+\mu^2)^2}dy \\
 &=\frac{9}{4}\frac{ (\mu^2+1) \ln(\mu^2)-2(\mu^2-1)}{\mu^2-1}, 
 \end{split}
\end{equation}

and 
$$
\int_{\Si^2}  u_{\frac{2}{3}, |\vec{a}|} d\omega = 
\frac{1}{4\pi}  \int_{\R^2}  \frac{3}{2} \ln (\frac{1+|y|^2 }{\mu^2+ |y|^2}) \frac{4}{(1+|y|^2)^2} dy +\ln(\mu^2) +\frac{1}{2} \ln  (\frac{2}{1+\mu^2} )
$$
$$= -\frac{3}{2} \frac{\mu^2 \ln(\mu^2) +\mu^2-1 }{\mu^2-1}+ \ln(\mu^2) +\frac{1}{2} \ln  (\frac{2}{1+\mu^2}).
$$
Then we can calculate  
$$
\min_{u \in \calM_{\vec{a}} } I_{\frac{2}{3}} (u)=I_{\frac{2}{3}} (u_{\frac{2}{3}, |\vec{a}|}) =0.
$$
Furthermore, 
\begin{equation}\label{upper}
\begin{split}
\min_{u \in \calM_{\vec{a}} } I_{\alpha} (u) &\le  I_{\alpha} (u_{\frac{2}{3}, |\vec{a}|} ) =(\alpha-\frac{2}{3}) \int_{\Si^2} |\nabla u_{\frac{2}{3}, |\vec{a}|}|^2 d\omega\\
&=(\alpha-\frac{2}{3})  \frac{9}{4|\vec{a}|}\bigl( -2 |\vec{a}|+  \ln \frac{1+|\vec{a}|}{1-|\vec{a}|} \bigr) .
\end{split}
\end{equation}

In particular,  if $\alpha < \frac{2}{3}$  we have that $  \min_{u \in \calM_{\vec{a}} } I_{\alpha} (u)  \to -\infty $  as $|\vec{a}| \to 1$.
This establishes the proof of Theorem \ref{main}.  $\blacksquare$
\vskip .2in


\section{Uniqueness and symmetry} 

For a better understanding of $ I_{\alpha} (u)$,   particularly for  $ \alpha \not =  \frac{2}{3}$,  we need to 
consider the minimizer $u_{\alpha, \vec{a}} $ in \eqref{minimizer}  more closely.   

First, we can rotate the coordinates properly so that $\beta_1=\beta_2=0$.   From \eqref{condition} we see that $\beta_3 \not =0$ if $\vec{a}  \not = (0, 0, 0)$.  Without loss of generality,  we assume  that $a_3 \ge 0$.
In view of \eqref{condition},  we have in particular $\rho= 1+ \beta_3 a_3 $ and  \eqref{minimizer} becomes

\begin{equation}\label{reduced}
\alpha \Delta u+e^{2u}(1+\beta_3 (a_3-x_3) )=1, \quad x \in \Si^2.
 \end{equation}
Also  \eqref{condition} is reduced to 
\begin{equation}\label{beta_3}
2 (\frac{1}{\alpha} -\frac{3}{2})  \beta_3 \int_{\Si^2} (x_3)^2 e^{2u_{\alpha, \vec{a}} } d\omega =
 2( \frac{1}{\alpha}-1)\rho a_3-\beta_3.
 \end{equation}
 
Using \eqref{beta_3} and the fact that
$$
(a_3)^2 =(\int_{\Si^2} x_3 e^{2u_{\alpha, \vec{a}}} )^2 d\omega \le \int_{\Si^2} (x_3)^2 e^{2u_{\alpha, \vec{a}} } d\omega \le 1,
$$
we can obtain
\begin{equation}\label{beta-small}
 \frac{a_3}{1-a_3^2} \le  \beta_3 \le  \frac{ 2(\frac{1}{\alpha}-1) a_3}{1-a_3^2}, \quad \hbox{ if }  \alpha \in (\frac{1}{2},  \frac{2}{3}]
 \end{equation}
 and
 \begin{equation}\label{beta-large}
 \frac{a_3}{1-a_3^2} \ge  \beta_3 \ge  \frac{ 2(\frac{1}{\alpha}-1) a_3}{1-a_3^2}, \quad \hbox{ if }  \alpha \in [\frac{2}{3}, 1].
 \end{equation}

 We first show that  when $\alpha \in (1/2, 1)$ is fixed,  the solution set  $\{u_{\alpha,  \vec{a}}\}  $ with parameters  $  \vec{a}$  has a simple structure  near  the trivial solution $u=0$ with $(\alpha,  \vec{a})=(\alpha, (0, 0, 0) )$.
 
  \begin{proposition}\label{uniqueness_a}
 Fix $\alpha \in (1/2, 1)$.   There is a  constant $\delta(\alpha)>0$ sufficiently small such that  when $ 0<|\vec{a}|<\delta(\alpha)$,  \eqref{reduced}  has  a unique solution  $u_{\alpha,  \vec{a}}$ in $\calM_{a}:= \calM_{\vec{a}}$,  which is therefore axially symmetric around $\vec{a}$.   
 \end{proposition}
 
 {\bf Proof. }   We only need to consider the special case $\vec{a}=(0,0, a_3) $ after a proper rotation.   From \eqref{beta-small} and \eqref{beta-large}, we know that 
 $
\beta \to 0,  \,\, \hbox{as} \,\, a_3 \to 0
$
and hence $u_{\alpha,  \vec{a}} $ converges to  the trivial solution $u=0$ as $a_3 $ goes to zero. 
Furthermore,  from \eqref{beta_3} we obtain
$$
\lim_{a_3 \to 0} \frac{\beta}{a_3}=3 (1-\alpha).
$$
Suppose there is a sequence of $\{a^{(k)} :=a^{(k)}_3,  k=1, 2, \cdots \}$ with $|a^{(k)}| \to 0$ as $k \to \infty$  such that 
$
\min_{u \in \calM_{\vec{a}} }  I_{\alpha} (u)
$
has two distinct solutions $u_1^{(k)}, u_2^{(k)}$  which satisfies \eqref{reduced} with $\beta_1^{(k)}, \beta_2^{(k)} $ respectively.

From \eqref{beta-small} and \eqref{beta-large}, we know that 
 $
\beta_1^{(k)}, \beta_2^{(k)}  \to 0,  \,\, \hbox{as} \,\, k  \to \infty
$
and hence $u_1^{(k)}, u_2^{(k)}$   converge  in $C^2 (\Si^2)$ to  the trivial solution $u=0$ as $k $ goes to infinity,  due to the uniqueness  result Theorem \ref{MAOTheorem}.
Furthermore,  from \eqref{beta_3} we obtain
$$
\lim_{k \to \infty} \frac{\beta_i^{(k)} }{a^{(k)} }=3 (1-\alpha), \quad i=1, 2.
$$

Let $ w_k:=u_1^{(k)}- u_2^{(k)}$,  and $m_k:=\|u_1^{(k)}- u_2^{(k)}\|_{L^{\infty}(\Si^2)}$. 
Note that   $w_k$ satisfies
$$
\alpha \Delta w_k+( e^{2u_1^{(k)}}-e^{2u_2^{(k)} }) \bigl(1+\beta_1^{(k)} (a^{(k)}-x_3) \bigr) = e^{2u_2^{(k)} } (\beta_2^{(k)}-\beta_1^{(k)} ) (a^{(k)}-x_3).
$$
Multiplying the above equation by $a^{(k)}-x_3$ and integrating on $\Si^2$, we can obtain
$$
|\beta_1^{(k)}-\beta_2^{(k)}| \le C \|u_1^{(k)}- u_2^{(k)}\|_{L^{\infty}(\Si^2)}
 $$
 for some positive constant $C$.  
 
 Hence, after taking a proper subsequence,  we know that  $w_k/m_k$ converges in $C^2(\Si^2)$  to  some nontrivial function $\phi \in C^2(\Si^2)$, 
 and $ (\beta_1^{(k)}-\beta_2^{(k)})/m_k$ converges to a constant $\mu_0 \in \R$.   Furthermore,  we have 
 $$
 \int_{\Si^2} \phi d\omega=  \int_{\Si^2} \phi x_3 d\omega=0 
 $$
 and 
 $$
 \alpha \Delta \phi + 2 \phi +\mu_0 x_3=0,  \quad x \in \Si^2.
 $$
 Multiplying the above equation by $x_3$ and integrating on $\Si^2$,  we obtain $\mu_0=0$. 
 Since $\alpha \in (1/2, 1)$ and the first  and second eigenvalue of the Laplacian on $\Si^2$  are $\lambda_1=2, \lambda_2=6$
 respectively,  this leads to  a contradiction.  The proposition is proven. $\blacksquare$
 
 Next we shall show a uniqueness result for $\alpha $ close to $\frac{2}{3}$ when $\vec{a}$ is fixed.  
 
 \begin{proposition}\label{uniqueness}
 Fix $\vec{a}\in B_1$.   There is a  constant $\delta(\vec{a}) >0$ sufficiently small such that  when $ |\alpha-\frac{2}{3}|<\delta $,  \eqref{reduced} 
has  a unique solution in $\calM_{a}:= \calM_{\vec{a}}$,  which is therefore axially symmetric.   
 \end{proposition}
 
 {\bf Proof. }  We only need to consider solutions to  \eqref{reduced}  in $\calM_{\vec{a}} $  with possible different $\beta_3$.   Assume the contrary that there is  a sequence of $\alpha_k,  k=1, 2, \cdots $ such that  $ \alpha_k \to \frac{2}{3} $ as $k \to \infty$,    and   \eqref{reduced}  has  distinct solutions $u_{1, k},  u_{2, k} $ corresponding to  possibly distinct values of   $\beta_3=\beta_k^1,  \beta_k^2$ and distinct rotations  $\vec{a}_k^1, \vec{a}_k^2$ of  $\vec{a} $ respectively.   It is well-known that these solutions are smooth and uniformly  bounded. 
By the uniqueness of solution to \eqref{eq-a},  it is  easy to see that $\vec{a}_k^1, \vec{a}_k^2$ converge to  $\vec{a} =(0, 0, a_3) $,  and  $ \beta_k^1,  \beta_k^2$ converge to $\beta_3=\frac{|\vec{a}|}{1-|\vec{a}|^2}=\frac{a_3}{1-a_3^2}$. 

  In view of \eqref{condition},  it is also easy to see that
 $$
 |\beta_k^1- \beta_k^2| \le C(|\vec{a}| )  |\alpha-\frac{2}{3}| \times \|u_{1, k}-u_{2, k} \|_{L^\infty (\Si^2)}
 $$
 for some positive constant $C$ depending only  on $|\vec{a}|  \in (0, 1)$. 
  
Let 
 $$\phi_k= \frac{u_{1, k}-u_{2, k}}{ \|u_{1, k}-u_{2, k} \|_{L^\infty (\Si^2)} }.
 $$
 It is standard to verify that $u_{1, k},  u_{2, k} $ converges to  $u_{\frac{2}{3}, \vec{a}}$   in $C^2(\Si^2)$ with 
 $\vec{a}=(0, 0, a_3)$,  as $k \to \infty$,  and  $\phi_k $ converges, after passing to a subsequence,   in $C^2(\Si^2)$ to $ \phi $ with $ ||\phi ||_{L^\infty (\Si^2)}=1$.  
 Furthermore, $\phi$ satisfies \eqref{constraint}  and the linearized equation
 \begin{equation}\label{linear3.5}
 \frac{2}{3} \Delta \phi +\frac{2(1-a_3^2)} {(1-a_3x_3)^2} \phi=0, \quad x\in \Si^2.
 \end{equation}
 Now consider the eigenvalue problem
  \begin{equation}\label{eigenvalue}
 \ \Delta \phi +\frac{\lambda (1-a^2)} {(1-a x_3)^2} \phi=0, \quad x\in \Si^2.
 \end{equation}
 for a fixed  $a \in (0, 1)$.
 Note that  for any $a \in (0, 1)$,   the transformation $ T_{a}:  \Si^2 \to \Si^2$ given by
 $$
  T_a(x)=\bigl(\frac{\sqrt{1-a^2} x_1}{ 1-a x_3},  \frac{\sqrt{1-a^2} x_1}{ 1-a  x_3} , \frac{x_3-a}{ 1-a  x_3}  \bigr)
 $$
 is a conformal transformation.  Indeed,  
 $$
 T_a= \Pi^{-1} \bigl( \sqrt{\frac{1-a}{1+a}}  \Pi ):   \Si^2 \to \Si^2
 $$
 and
 $$
 det (d T_a)=  \frac{(1-a^2)} {(1-a x_3)^2}.
 $$
 Then we  observe  that $P(T_a(x))$ is an eigenfunction to \eqref{eigenvalue} if and only if $P(x)$ is a spherical harmonics.  Therefore,  \eqref{eigenvalue} has only eigenvalues  $\lambda=m(m+1)$ for  a  nonnegative integer $m$.   This leads to a contradiction to \eqref{linear3.5} since $\lambda=3$ is not an eigenvalue of \eqref{eigenvalue}.  The proof is complete. $\blacksquare$
 
 We can show the axial symmetry of a minimizer to  $\min_{u \in \calM_{\vec{a}} }  I_{\alpha} (u)$ for most  cases,  though it is still not completely resolved whether a given minimizer is always axially symmetric.
 
  \begin{proposition}\label{symmetry}
 Fix $\vec{a}\in B_1$,   assume that  for $\alpha >1/2$,  $u_{\alpha, \vec{a}} $   is a solution to \eqref{reduced} in $\calM_{a}:= \calM_{\vec{a}}$. 
Then $u_{\alpha, \vec{a}} $ must be axially symmetric  when either  i)  $\alpha \in (1/2, 2/3] $ or  ii) $\alpha \ge 1 $ or  iii) $ \alpha \in  (2/3, 1) $ and $
|\vec{a} |  \le \frac{1-\alpha}{2\alpha-1}.  $
 \end{proposition}

 {\bf Proof. } 
For this purpose,  we choose the stereographic project $ \Pi : \Si^2 \to \R^2$ from the north pole $N=(0, 0, 1)$.
By  \eqref{beta-small} and \eqref{beta-large},  we have
$$
0\le \beta_3 < \frac{1}{1-a_3},  \quad  \rho- \beta_3>0.
$$

 Set 
\[w_{\alpha, \vec{a}} (y):=u_{\alpha, \vec{a}} (\Pi^{-1}(y))-\frac{1}{\alpha}  \ln(1+|y|^2 )  +\frac{1}{2} \ln( \frac{4(\rho-\beta_3)} {\alpha}) \ \ \hbox{for}\ \ y\in \R^2.\]

Let $\mu$ be a positive constant with $\mu^2=\frac{\rho+\beta_3}{\rho-\beta_3}>1 $.    Then $w_{\alpha, \vec{a}} $ satisfies 
\begin{equation}\label{general-plane}
\Delta w+ \calK(|y|) e^{2w}=0 \ \ \hbox{in} \ \ \R^2
\end{equation}
and 
\begin{equation}\label{general-total}
\frac{1}{2 \pi} \int_{\R^2}\calK(|y|) e^{2w} dy= \frac{2}{\alpha}
\end{equation}
where
$$
\calK(|y|):= (\mu^2+|y|^2) (1+|y|^2)^{\frac{2}{\alpha}-3}.
$$
i) When $\frac{1}{2} <\alpha <\frac{2}{3}$,   it is easy to see that  $\calK(|y|) $ satisfies $(K1)-(K2)$ with $l=\frac{2}{\alpha}-2$. 
By Proposition \ref{th-general},   we know that $w_{\alpha, \vec{a}} (y)$ is radially symmetric and hence 
$u_{\alpha, \vec{a}} (y)$ must be axially symmetric and  $a_1=a_2=0$.

ii) When $\alpha >  1$,  then $\calK >0$ is not constant and decreasing in $r=|y|$.  The standard moving plane method can lead to the radial symmetry of $w_{\alpha, \vec{a}} (y)$.  Indeed,  the  radially symmetric solution is also  unique (see  Theorem 1.4  of \cite{Lin1-MR1770683}).  

 When $\alpha=1$,  by \eqref{beta_3},  we know that $\beta_3=0$ and  hence \eqref{reduced} becomes \eqref{standardPDE} with 
$\alpha=1$.   It is well known that there is a unique solution  to $u_{1, \vec{a}} \in \calM_{\vec{a}}$  which is axially symmetric about $\vec{a}$.  

iii) When $\frac{2}{3} <\alpha <1$,  if $ \mu^2 (3-\frac{2}{\alpha}) \le 1$ we have 
$$
\Delta \ln \calK(|y|)= \frac{ 4 [\mu (r^2+1)+\sqrt{3-\frac{2}{\alpha}} (r^2+\mu^2) ]}  { (r^2+1)^2 (r^2+\mu^2)^2} \Bigl(\mu (r^2+1)-\sqrt{3-\frac{2}{\alpha}} (r^2+\mu^2)\Bigr) \ge 0
$$
 and hence  $ (K1)-(K2) $ are   satisfied. 
In particular,  in view of \eqref{beta-large},  $(K1)-(K2)$ hold when 
\begin{equation}\label{a_3}
a_3 =|\vec{a}| \le \frac{1-\alpha}{2\alpha-1}.
\end{equation} 

By Proposition \ref{th-general} with $l= \frac{2}{\alpha}-2$ ,  under the condition \eqref{a_3},  $w_{\alpha, \vec{a}} (y)$ must be  radially symmetric   and hence  $u_{\alpha, \vec{a}} (x)$ must be axially symmetric, $a_1=a_2=0$. $\blacksquare$
\vskip .1in

\section{Estimates of the minimum of $I_{\alpha} $ on $\calM_{\vec{a}}$} 
In this section,  we shall estimate for $\alpha \in (1/2, 1)$ 
\begin{equation}\label{minimizer_a}
m(\alpha, a):= \inf_{ u \in \calM_{\vec{a}}, |\vec{a}|=a}   I_{\alpha} ( u ),
\quad \forall a \in [0, 1).
\end{equation}

In view of Proposition \ref{compact},  we know that  $ m(\alpha, a) $  is a continuous function of $a \in [0, 1)$ for any fixed $\alpha \in (1/2, 1)$.

We have the following estimates.

\begin{theorem}\label{energy}
There hold pointwise in $a \in [0, 1)$
\begin{equation} \label{lowerbound} 
  m(\alpha, a) \ge \left\{
\begin{aligned}
   &(\frac{2}{\alpha}-3) \ln(1-a^2) , \quad \alpha \in (1/2, 2/3),   \\
  &\alpha  (\frac{1}{\alpha}-\frac{3}{2}) \ln(1-a^2), \quad \alpha \in (2/3, 1).
\end{aligned}
\right.
\end{equation}
and  
\begin{equation}\label{upperbound}
  m(\alpha, a) \le \left\{
\begin{aligned}
 &(\frac{2}{\alpha}-3) \ln(1-a^2), \quad \alpha \in (2/3, 1), \\
 & \frac{3\alpha}{2a} (\frac{1}{\alpha}-\frac{3}{2}) \bigl( \ln(1-a^2)-2(\ln(1+a)- a) \bigr), \quad \forall \alpha \in (1/2, 1).  
\end{aligned}
\right.
\end{equation}

There also holds asymptotically as $a \to 1$
\begin{equation}\label{upper-approx}
    m(\alpha, a) \le  (\frac{1}{\alpha}- \frac{3}{2}) \ln(1-a^2)  \bigl( 1 + o(1) \bigr),\,\, \quad  \alpha \in (1/2, 1).
  \end{equation}
 \end{theorem}

{\bf Proof. }

We first recall Onofri's inequality 
 \begin{equation}\label{onofri}
 F_1(u)=\int_{\Si^2}|\nabla u| ^2 d\omega+ 2\int_{\Si^2}u d\omega-\log \int_{\Si^2}e^{2u}d \omega \ge 0, \quad u \in H^{1} (\Si^2).
 \end{equation}
 
In view of \eqref{ineq} and \eqref{onofri},  it is easy to see by interpolation  that   for  $\alpha \in (2/3, 1]$
 \begin{equation}\label{lower-large}
  I_{\alpha} (u ) \ge \alpha  (\frac{1}{\alpha}-\frac{3}{2}) \ln(1-a^2), \quad \forall u \in \calM.
\end{equation}

As we know from previous discussion,  there is  a minimizer  $u_{\alpha, a} $ to the minimization problem \eqref{minimizer_a},  
which is a solution to \eqref{reduced} with $a_3=a,  \beta_3=\beta_3(a)$ satisfying \eqref{beta-small} and \eqref{beta-large}.
Also from  Proposition \ref{uniqueness_a},   the solution  $u_{\alpha, a} $ forms a curve smooth curve parametrized by $a \in (0, \delta(\alpha))$. 
Furthermore, the linearized operator of \eqref{reduced} is a Fredholm operator on the tangent space of $\calM_{\vec{a}}$ at any solution $u$ of \eqref{reduced} on $\calM_{\vec{a}}$.  
By the compactness of solutions of  \eqref{reduced}  for $ a_3=a \in [0, 1-\epsilon]$  for any fix $\epsilon \in (0, 1)$ and the analyticity of 
equation   \eqref{reduced} in term of  $u$,   it can be  shown by the global bifurcation theory (see, e.g, Theorem 9.1.1 in \cite{MR1956130}) that any solution set  of \eqref{reduced} can be extended globally with either 
$a \to 0,  1$ or being  a closed loop.  In particular, by Proposition \ref{uniqueness_a} there exists a branch of solution set which extends to $a=0$ in one direction  and to $a=1$ in the other direction. Note that we do not know in general the uniqueness of the solution for a fixed $a \in (0, 1)$,  there might be more branches,  and each branch of solutions might contain portions which are not minimizers of \eqref{minimizer_a}.  

Nevertheless, by the compactness result Proposition \ref{compact} again, there are only finite numbers of smooth branches of solutions to \eqref{reduced} for  $a \in [0, 1-\epsilon] $ with $\epsilon>0$,  and we can find  a piecewise smooth solution curve  $u_{\alpha, a(\tau)}, \tau \in (0, \infty) $  to   \eqref{reduced} in $\calM_{a}:= \calM_{\vec{a}}$  with $\vec{a}=(0, 0, a)$  such that  $ a(\tau) \to 1$ as $\tau \to \infty$.  Furthermore,   the singular set 
$$
\calS:= \{ \tau:  a'(\tau)\,  \hbox{ does not exists ;  or  } \, a'(\tau) =0  \, \hbox{;  or} \,   \, \frac{\partial u_{\alpha, a(\tau)}}{\partial a}  \hbox{ does not exists in }  \, C^2 (\Si^2)\}
$$
 does not have accumulative point. 
  For a fixed $a \in(0, 1)$, there is  some $\calT>0$ depending on $\alpha$ such that  $a(\calT) =a$. 
  We have
  $$
  a'(\tau) \int_{\Si^2} e^{ 2 u_{\alpha, a(\tau)}}  \frac{\partial u_{\alpha, a(\tau)}}{\partial a}  d\omega  =0, \quad \tau \in [0, \calT]\setminus \calS
$$
  and
  $$
  a'(\tau) \int_{\Si^2} e^{ 2 u_{\alpha, a(\tau)}}  x_3  \frac{\partial u_{\alpha, a(\tau)}}{\partial a}  d\omega  =a'(\tau), \quad \tau \in [0, \calT]\setminus \calS.
  $$
  
   Now  using \eqref{reduced} and the above equalities we obtain 
 $$
\frac{\partial   I_{\alpha} (u_{\alpha, a(\tau)} ) }{\partial \tau}= -2 (\beta(\tau) -\frac{a(\tau)}{1-a^2(\tau)} )  a'(\tau),  \quad \tau \in [0, \calT]\setminus \calS.
$$
Using similar arguments,  we can find a  solution  curve of \eqref{minimizer_a},  still denoted by  $u_{\alpha, a(\tau)},  \tau \in [0, \calT], $ which
is piecewise smooth,  but may have finite discontinuous points $\tau_i,  i=1, 2, \cdots,  N$ with $\tau_0=0, \tau_{N}= a$.   Moreover,  it can be chosen that $a(\tau) $ is continuous,  and $u_{\alpha, a(\tau)} $ has both left limit and right limit at $\tau_i,   i=0, 1, 2, \cdots, N$ in $C^2(\Si^2)$ and $ I_{\alpha} (u_{\alpha, a(\tau)})=m(\alpha, a(\tau) )$ is continuous.  Furthermore,  $a'(\tau) >0,  \tau \in (\tau_i, \tau_{i+1}),  i=0, 1, 2, \cdots, N$.

Hence,  using \eqref{beta-small} and \eqref{beta-large} with $a_3=a(\tau) $ and $\beta_3=\beta_3(\tau)$,  in view of $ I_{\alpha} (u_{\alpha, 0} ) =0$  we have for $ \alpha \in (1/2, 2/3]$
$$
0 \ge I_{\alpha} (u_{\alpha, a} ) =\int_{0}^{\calT} \frac{\partial   I_{\alpha} (u_{\alpha, a(\tau)} ) }{\partial \tau} d\tau
$$
$$
=\sum_{i=1}^{N}\int_{\tau_{i-1}}^{\tau_i} \frac{\partial   I_{\alpha} (u_{\alpha, a(\tau)} ) }{\partial \tau} d\tau  \ge \sum_{i=1}^{N}\int_{\tau_{i-1}}^{\tau_i}  (\frac{2}{\alpha}-3) \frac{2 a(\tau) a'(\tau)}{1-a^2(\tau)}  d \tau
$$
 and for  $\alpha \in (2/3, 1]$
$$
0 \le  I_{\alpha} (u_{\alpha, a} ) =\int_{0}^{\calT} \frac{\partial   I_{\alpha} (u_{\alpha, a(\tau)} ) }{\partial \tau} d\tau
$$
$$= \sum_{i=1}^{N}\int_{\tau_{i-1}}^{\tau_i} \frac{\partial   I_{\alpha} (u_{\alpha, a(\tau)} ) }{\partial \tau} d\tau \le \sum_{i=1}^{N}\int_{\tau_{i-1}}^{\tau_i}  (\frac{2}{\alpha}-3) \frac{2 a(\tau) a'(\tau)}{1-a^2(\tau)}  d \tau.
$$
Hence the first inequalities in both  \eqref{lowerbound} and \eqref{upperbound} are proven.

 Next, we will estimate  $m(\alpha, a(\tau))$ from above by using suitable auxiliary functions.
Define
 \begin{equation}\label{test}
{\tilde u}_{\alpha, \mu}(x)={\tilde u}_{\alpha, \mu}(\Pi^{-1}(y)): =\frac{1}{\alpha} \ln\frac{1+|y|^2 }{\mu^2+ |y|^2}, \quad y \in \R^2.
\end{equation}

Direct computations show that 
$$
\int_{\Si^2} |\nabla {\tilde u}_{\alpha, \mu}  |^2 d\omega = 
\frac{1}{4\pi}  \int_{\R^2} |\nabla {\tilde u}_{\alpha, \mu}(\Pi^{-1}(y))|^2dy = \frac{1}{\alpha^2 (\mu^2-1)} \bigl( 2(1-\mu^2)+ (\mu^2+1) \ln(\mu^2)\bigr), 
$$
$$
\int_{\Si^2} {\tilde u}_{\alpha, \mu}  d\omega = 
\frac{1}{4\pi}  \int_{\R^2} {\tilde u}_{\alpha, \mu}(\Pi^{-1}(y)) \frac{4}{(1+|y|^2)^2} dy = -\frac{1}{\alpha  (\mu^2-1)} \bigl( (1-\mu^2)+ \mu^2\ln(\mu^2)\bigr), 
$$
and

$$
\int_{\Si^2} e^{2{\tilde u}_{\alpha, \mu} }  d\omega = \frac{1}{4\pi}  \int_{\R^2}e^{2 {\tilde u}_{\alpha, \mu}(\Pi^{-1}(y)) } \frac{4}{(1+|y|^2)^2} dy=
\frac{1-\mu^{2-4/\alpha}} { ( \frac{2}{\alpha}-1)(\mu^2-1)}
$$

$$
\int_{\Si^2} e^{2{\tilde u}_{\alpha, \mu} } x_3 d\omega = \frac{1}{4\pi}  \int_{\R^2}e^{2 {\tilde u}_{\alpha, \mu}(\Pi^{-1}(y)) } \frac{4(|y|^2-1)}{(1+|y|^2)^3} dy
$$
$$
=\frac{ ( \frac{2}{\alpha}) \mu^{4-4/\alpha}  +(  \frac{2}{\alpha}-2) (\mu^2-\mu^{2-4/\alpha}) - \frac{2}{\alpha} } { ( \frac{2}{\alpha}-1) ( \frac{2}{\alpha}-2)(\mu^2-1)^2}
$$

Hence,   we have 
\begin{equation}
 I_\alpha({\tilde u}_{\alpha, \mu})= (\frac{3}{2} - \frac{1}{\alpha} ) \ln( \mu^2 ) \bigl( 1 + o(1) \bigr) \end{equation}
 as $a \to 1 $  and $\mu \to \infty$.
 
 Now we compute the center of mass of ${\tilde u}_{\alpha, \mu}$ 
 $$
a_{\alpha, \mu}:= \frac{\int_{\Si^2} e^{2{\tilde u}_{\alpha, \mu} } x_3 d\omega}{\int_{\Si^2} e^{2{\tilde u}_{\alpha, \mu} }  d\omega}=1-\frac{2(\frac{2}{\alpha}-1)(1-\mu^2) +2(\mu^{\frac{4}{\alpha}-2}-1)}{(\frac{2}{\alpha} -2) (\mu^2 -1) (\mu^{\frac{4}{\alpha}-2}-1) }.
$$
It is easy to see that  for $\alpha \in (1/2, 1)$
$$
a_{\alpha, \mu} \to 1 \,\, \hbox{as }\,\,  \mu \to \infty;  \quad a_{\alpha, \mu} \to 0 \,\, \hbox{as }\,\,  \mu \to 1.
$$
Hence  for any fixed $a \in (0, 1)$ there is at least a  positive number  $\mu(a) $ such that  $a_{\alpha, \mu(a)} =a$ and
$$
\mu^2(a) =   \frac{\alpha}{(1-\alpha)(1-a)}  \bigl( 1 + o(1) \bigr),  \,\, \hbox{as }\,\,  a \to 1. 
$$
Letting $\mu=\mu(a)$,  we can choose a constant $c=c(a)$ such that $v_{\alpha, a} := u_{\alpha, \mu(a)}+c(a) \in \calM_{a}$ and
$$
I_\alpha(v_{\alpha, a} ) =  (\frac{1}{\alpha}- \frac{3}{2}) \ln(1-a^2)  \bigl( 1 + o(1) \bigr) \,\, \hbox{as }\,\,  a \to 1. 
$$
 Hence  \eqref{upper-approx} holds.

 Therefore,  the proof is complete.   $\blacksquare$

 \begin{remark}
 In view of \eqref{onofri},   we can  then derive   a lower bound of  the gradient  $L^2$ norm of   a minimizer  $u_{\alpha, a} $ of \eqref{minimizer_a}  asymptotically
 $$
 \int_{\Si^2}|\nabla u_{\alpha, a}| ^2 d\omega \ge - \frac{\ln(1-a^2)  }{\alpha} \bigl( 1 + o(1) \bigr)\,\, \hbox{as }\,\,  a \to 1;
 $$
 While for $\alpha \in (1/2, 2/3)$,  using \eqref{ineq1} a better lower bound can also be obtained
 $$
 \int_{\Si^2}|\nabla u_{\alpha, a}| ^2 d\omega \ge - \frac{3 \ln(1-a^2)  }{2 \alpha} \bigl( 1 + o(1) \bigr)\,\, \hbox{as }\,\,  a \to 1. 
 $$
 By \eqref{ineq} and the above gradient estimates,  we  can obtain  for $\alpha \in (2/3, 1)$  an  asymptotical lower bound 
 $$
 m(\alpha, a) \ge  (\alpha- \frac{2}{3})  \int_{\Si^2}|\nabla u_{\alpha, a}| ^2 d\omega \ge  \frac{2}{3} (\frac{1}{\alpha}-\frac{3}{2}) \ln(1-a^2)  \bigl( 1 + o(1) \bigr)\,\, \hbox{as }\,\,  a \to 1. 
 $$
 However,  this  energy lower bound is not as good as \eqref{lower-large}. 
  When $\alpha \in (1/2, 2/3)$,  
 similar  energy lower bound does not follow immediately from the gradient lower bound.  If it did,  it would be a sharp one as it  would coincide with the upper bound.   Nevertheless,  it is expected that 
$$
 m(\alpha, a)=(\frac{1}{\alpha}-\frac{3}{2}) \ln(1-a^2)  \bigl( 1 + o(1) \bigr)\,\, \hbox{as }\,\,  a \to 1. 
 $$

Similarly using \eqref{onofri}  and \eqref{upperbound},  we can obtain,   when $\alpha \in (2/3, 4/5)$,  a lower bound of the gradient  $L^2$ norm of   a minimizer  $u_{\alpha, a} $ of \eqref{minimizer_a}  point wisely in $a$
  $$
 \int_{\Si^2}|\nabla u_{\alpha, a}| ^2 d\omega \ge  - \frac{4-5\alpha}{2\alpha (1-\alpha)}  \ln(1-a^2), \quad  \,\, \forall a \in [0, 1], 
 $$
  which leads to the following  lower bound point wisely in $a$
 $$
 m(\alpha, a) \ge  (\alpha- \frac{2}{3})  \int_{\Si^2}|\nabla u_{\alpha, a}| ^2 d\omega \ge  -   \frac{(4-5\alpha)(2\alpha-3)}{4\alpha (1-\alpha)}  \ln(1-a^2).
$$
We note that the upper bound  in \eqref{upperbound} is not optimal  which  leads  to  the technical condition  $\alpha<4/5$ in the above estimates  instead of the  natural range up to $\alpha<1$.

 On the other hand, using  \eqref{ineq1}  and \eqref{upperbound},  we can also derive an upper bound of the gradient $L^2$ norm of 
 a minimizer  $u_{\alpha, a} $ of \eqref{minimizer_a}    point wisely in $a$  when $\alpha \in (2/3, 1)$
 $$
 \int_{\Si^2}|\nabla u_{\alpha, a}| ^2 d\omega \le  - \frac{3}{\alpha}  \ln(1-a^2), \quad  \,\, \forall a \in [0, 1]. 
 $$
 Similarly, using \eqref{ineq1} and \eqref{upper-approx},  we can  derive an upper bound of the gradient $L^2$ norm of a minimizer  $u_{\alpha, a} $ of \eqref{minimizer_a}  asymptotically in $a$
  when $\alpha \in (2/3, 1)$
  $$
 \int_{\Si^2}|\nabla u_{\alpha, a}| ^2 d\omega \le  - \frac{3}{2\alpha}  \ln(1-a^2)  \bigl( 1 + o(1) \bigr)\quad  \hbox{as }\,\,  a \to 1. 
 $$
 Similar upper bounds seem not follow immediately when $\alpha \in (1/2, 2/3)$.

 \end{remark}

\begin{remark}
The following technical questions still remain open:

1)  Should $u_{\alpha, \vec{a}} (x)$ always be axially symmetric for all $\alpha \in (\frac{2}{3}, 1)$ and $\vec{a} \in B_1$?  

2)  Is the minimizer  $u_{\alpha, \vec{a}} (x)$ unique determined?   In particular, is $\beta$ uniquely determined?   We know that if $\beta$ is uniquely determined by $\alpha $ and $\vec{a}$, then the axially symmetric solution $u_{\alpha, \vec{a}} (y) $  is unique. 

3)  Fixed $\alpha \in (\frac{1}{2}, 1), \vec{a} \in B_1$,  for any given
 $\vec{\beta}=\beta_3 \vec{a} /|\vec{a}|,  0 <\beta_3<\frac{1}{1- |\vec{a}|}, \rho= 1+ \beta_3 |\vec{a}| $,  there is a unique  axially symmetric solution $u$ to \eqref{minimizer} with the corresponding 
 $w$ solving \eqref{general-plane}- \eqref{general-total}, following Theorem 1.5 of \cite{Lin1-MR1770683}.   However, it is not clear whether  
 the center of mass  $A \vec{a}/|\vec{a}|$ of $u$ is still $\vec{a}$ and the total mass $M$  is still $1$ or not.   Certainly for some such $\beta_3, \rho$,  we should have $A \not =|\vec{a}|, M \not =1$ since otherwise \eqref{beta-small} or \eqref{beta-large} should hold.  This implies that a solution to  \eqref{general-plane}-\eqref{general-total} may not   a solution to  the minimizing problem  $
\min_{u \in \calM_{\vec{a}} }  I_{\alpha} (u)
$.
  Nevertheless,  we note that  $M [1+\beta_3(|\vec{a}|-|A| )]=1$ still holds.

4)  Can we compute or estimate more accurately  $ m(\alpha, a).$
We need to get more information on the minimizer.  In particular,  what is the asymptotic behavior of the solution $ u_{\alpha, \vec{a}} $ as $ \vec{a} \in B_1$ goes to the unit sphere?  If we have
a detailed profile of the solution in different regions, we might get an optimal asymptotic estimate of  $ m(\alpha, a).$

We note that the answers to  all above questions are known for $\alpha =\frac{2}{3}$ as shown in 
Theorem \ref{main}.  The technical questions and  other related problems will be studied in a forthcoming paper. 

\end{remark}

\vskip .2in

\section{Second Variation of $I_{\alpha}$} 

Now we consider  another technical aspect of  $I_\alpha$:  the second variation of $I_\alpha$ in $H^1(\Si^2)$, in an effort to understand
$I_{\alpha}$ better.

 Fixed a solution $u \in \calH$  
to  \eqref{simple} and \eqref{a_i},  for any $\phi \in H^1(\Si^2)$  we have 
$$
D^2 I_{\alpha} (u) (\phi, \phi)=
  2 \alpha \int_{\Si^2} |\nabla \phi|^2 d\omega + \frac{ 8} {(1-|\vec{a}|^2 )^2}
  \bigl(\int_{\Si^2} e^{2u} (1-\vec{a} \cdot x ) \phi  d\omega \bigr)^2 
   $$
$$
- \frac{4  } {1-|\vec{a}|^2 } \bigl( \int_{\Si^2} e^{2u}  (1-\vec{a} \cdot x )\phi^2   d\omega+(\int_{\Si^2}  e^{2u} \phi  d\omega)^2 -\sum_{i=1}^{3} (\int_{\Si^2}  e^{2u}  x_i   \phi  d\omega)^2 \bigr).
$$

In particular, at $u \equiv 0$  for any $\phi \in H^1(\Si^2)$ 
$$
D^2 I_{\alpha} (0) (\phi, \phi)=2 \alpha \int_{\Si^2} |\nabla \phi|^2 d\omega -4 \int_{\Si^2}  \phi^2  d\omega +4 ( \int_{\Si^2}  \phi  d\omega)^2 +4 \sum_1^3(\int_{\Si^2}   x_i\phi d\omega)^2.
$$
Let $\phi =\sum_{n=0}^{\infty} b_n \phi_n  $  where $\{\phi_n\}_0^\infty  $ is an orthonormal basis of $ H^1(\Si^2)$ formed by spherical harmonics in an increasing order of eigenvalues $\lambda_n$.   Note that
$\phi_0 =1, \phi_i=  \sqrt3 x_i,  i=1, 2, 3 $ and $\lambda_i=0, \lambda_1=\lambda_2=\lambda_3=2$ and $\lambda_4=6$. 

Then $  \forall \phi \in H^1(\Si^2)$
$$
D^2 I_{\alpha} (0) (\phi, \phi)= 2\alpha  \sum_{n=0}^\infty \lambda_n b_n^2 -4 \sum_{n=0}^\infty b_n^2 + 4 b_0^2  +\frac{4}{3}  \sum_{n=1}^{3} b_n^2  
$$
$$
\ge  (4\alpha -\frac{8}{3}) \sum_{n=1}^{3} b_n^2 + (12 \alpha -4) \sum_{n=4}^\infty b_n^2.
$$

In particular,  the  linearized equation of equation \eqref{simple} at $u \equiv 0$ is
\begin{equation}\label{linear}
\alpha \Delta \phi +2\phi - 2\int_{\Si^2}  \phi  d\omega- 2 \sum_{i=1}^{3} x_i  \int_{\Si^2}   x_i\phi d\omega=0  \quad \hbox{ on }  \quad \Si^2.
\end{equation}
It has a kernel $\calK_\alpha=\{ \phi_0\} $ if $\alpha \not  =\frac{2}{3}$ and  
$\calK_{\frac{2}{3}}=\{ \phi_i,  i=0, 1, 2, 3\} $ when $\alpha = \frac{2}{3}$.

Hence, it is easy to conclude from the above discussion  the following 
\begin{proposition}\label{second variation}

 i)  when $\alpha \in (\frac{2}{3}, 1)$,  $D^2 I_{\alpha} (0) (\phi, \phi)\ge 0$ 
and the equality holds only when $\phi $ is a constant function;  In particular, $D^2 I_{\alpha} (0)$
is positive definite when restricted to $\calH$; 
 
 ii) when $\alpha =\frac{2}{3} $,  $D^2 I_{\alpha} (0) (\phi, \phi)\ge 0$ 
and the equality holds only when $\phi $ is  expanded by  $\phi_0 =1, \phi_i=  \sqrt3 x_i,  i=1, 2, 3 $;

 iii) when $\alpha <\frac{2}{3} $,   $D^2 I_{\alpha} (0)$ is not non-negative.  In particular,  $I_\alpha (u) <0$ for some $u \in H^1(\Si^2).$. 
 
\end{proposition}
 This  fact gives a simple explanation of the critical value of $\alpha$ being $2/3$, compared to Theorem \ref{main}.

Finally, we shall look at the second variation of $I_{\alpha}$ at the nontrivial  explicit solution  when $\alpha=\frac{2}{3}$.
We can rewrite the solution as  
$$
u_{\frac{2}{3}, \vec{a}}(x)=-\frac{3}{2}\ln{(1-\vec{a} \cdot x)} +\ln (1-|\vec{a}|^2), \quad x \in \Si^2.
$$
Then, for any $\phi \in H^1(\Si^2),$   from Theorem \ref{main} we have 
$$
D^2 I_{ \frac{2}{3} } (u_{\frac{2}{3}, \vec{a}}) (\phi, \phi)
$$
$$
=\frac{4}{3} \int_{\Si^2} |\nabla \phi|^2 d\omega 
+8 (1-|\vec{a}|^2)^2 \bigl( \int_{\Si^2} \frac{ \phi}{(1-\vec{a} \cdot x)^2} d\omega \bigr)^2 - 4(1-|\vec{a}|^2)  \int_{\Si^2}  \frac{ \phi^2 }{(1-\vec{a} \cdot x)^2}   d\omega
$$
$$- 4(1-|\vec{a}|^2)^3  \bigl( (\int_{\Si^2}   \frac{\phi }{(1-\vec{a} \cdot x)^3}   d\omega)^2- \sum_{i=1}^3 (\int_{\Si^2}   \frac{ x_i\phi }{(1-\vec{a} \cdot x)^3} d\omega)^2 \bigr) \ge 0.
$$
In particular, if we consider the second variation of $I_{\frac{2}{3}} $ at $u_{\frac{2}{3}, \vec{a}} $ on
$\calM_{\frac{2}{3}, \vec{a}}$,  we only need to deal with  $\phi \in H^1(\Si^2)$ with
\begin{equation}\label{constraint}
\int_{\Si^2}   \frac{\phi }{(1-\vec{a} \cdot x)^3} d\omega=0, \quad \quad \int_{\Si^2}   \frac{ x_i\phi }{(1-\vec{a} \cdot x)^3} d\omega=0, \,\, i=1, 2, 3. 
\end{equation}

In this setting, it  also holds that 
$$
\int_{\Si^2}  \frac{ \phi }{(1-\vec{a} \cdot x)^2}   d\omega=0.
$$
Hence we have  
$$
D^2 I_{ \frac{2}{3} } (u_{\frac{2}{3}, \vec{a}}) (\phi, \phi)=\frac{4}{3} \int_{\Si^2} |\nabla \phi|^2 d\omega 
 - 4(1-|\vec{a}|^2)  \int_{\Si^2}  \frac{ \phi^2 }{(1-\vec{a} \cdot x)^2}   d\omega \ge 0.
$$

\vskip .2in

\section{Monotonicity}

In this section,  we shall discuss and  prove the first monotonicity formula of the analogue of the Szeg\"o Limit theorem on $\Si^2$.  
 Following \cite{GS-1958} (Section 2.1-2.2, Chapter 2), for any given function  $f\ge 0,  f \not \equiv 0$, we denote a measure 
 $ d\nu= f d\omega $ on $\Si^2$ and orthogonalize the functions $f_0=1, f_1=\sqrt{3} x_1, f_2=\sqrt{3}x_2, f_3=\sqrt{3}x_3$
 with respect to this measure.  Note that $f_0, f_1, f_2, f_3$ 
 form an orthonormal basis for spherical harmonics of order less than or equal to $1$ with respect to the measure $d\omega$.  

 Differing from the case on $S^1$ as discussed in \cite{GS-1958},  we only construct  $\phi_0, \phi_1, \phi_2, \phi_3$ such that  $\phi_0, \phi_1$ form an orthonormal basis for functions generated by $f_0, f_i$ respectively for $i=1, 2, 3$.  

 Denote the inner product on $L^2(\Si^2, d\nu)$ by $< , >$.
 We define
 $$ D_0= < f_0, f_0>=\int_{S^2} fd\omega $$
 and
$$
 D_{0, i} =det
 \begin{pmatrix}  
  < f_0, f_0> & < f_0, f_i>  \\
  < f_i, f_0> & < f_i, f_i>
  \end{pmatrix} , \quad i =1, 2, 3
  $$
  and 
  $$
  D_1= \frac{1}{3} \sum_{i=1}^3 D_{0, i}.
  $$

  It is easy to see that  $\phi_0= D_0^{-1/2} f_0$  and 
  $$
  \phi_i = (D_0D_{0, i})^{-1/2}\begin{vmatrix}  
  < f_0, f_0> & < f_0, f_i>  \\
   f_0 & f_i
  \end{vmatrix} = l_{i, 0} f_0+ l_{i, 1} f_i
  $$
  where
  $$
  l_{i, 0}= -(D_0D_{0, i})^{-1/2} < f_0, f_i>,\quad i=1, 2, 3,
  $$
  $$
  l_{i, 1}=(D_0 D_{0, i})^{-1/2} < f_0, f_0>=(\frac{D_0}{D_{0, i}})^{1/2}, 
  \quad i=1, 2, 3.
  $$

Now we state the following stage one monotonicity relation on $\Si^2$, which may be considered as the counterpart on $\Si^2$ of the Szeg\"o monotonicity theorem 
on $S^1$. 

\begin{proposition}\label{monotone}
We have 
$$
D_1= (\int_{\Si^2} f d\omega )^2- \sum_{i=1}^3 (\int_{\Si^2} f x_id\omega )^2 \ge0
$$
and 
\begin{equation}\label{monotonicity}
\ln D_0 -\int_{\Si^2} f d\omega  \le  \ln D_1 -2\int_{\Si^2} f d\omega.
\end{equation}

  \end{proposition}

{\bf Proof}  It is easy to see that 
$$
\sum_{i=1}^3 (\int_{\Si^2} f x_id\omega )^2 \le \sum_{i=1}^3 
\int_{\Si^2} f d\omega \cdot \int_{\Si^2} f x_i^2d\omega  \le 
(\int_{\Si^2} f d\omega)^2.
$$
To show \eqref{monotonicity}, we shall follow the proof of Theorem a of \cite{GS-1958} and prove
$$
\mu_i:= \inf_{a_i\in R} \int_{\Si^2} |a_i+ f_i|^2 d \nu 
$$ 
is attained at $a_i= \frac{l_{i,0}}{l_{i, 1}} $ and
$$
\mu_i=l_{i,1}^{-2} =\frac{D_{0, i}}{D_0}, \quad i=1, 2, 3.
$$
This can be seen from the folowing two facts.   First, 
$$
\mu_i \le \int_{\Si^2} |\frac{l_{i,0}}{l_{i, 1}} + f_i|^2 d \nu
\le l_{i,1}^{-2}  <\phi_i, \phi_i> = l_{i,1}^{-2} 
$$
Second,   if we write  $a_i+f_i= b_0 \phi_0 +b_i \phi_i$ where $b_i  l_{i, 1}=1$,  hence
$$
\int_{\Si^2} |a_i+ f_i|^2 d \nu = b_0^2+ b_i^2 \ge b_i^2= l_{i,1}^{-2}.
$$
Therefore
\begin{align*}
\frac{D_1}{D_0}= &\frac{1}{3} \sum_{i=1}^3 \mu_i \\
= &\frac{1}{3}  \inf_{a_i\in R} \int_{\Si^2} \sum_{i=1}^3 |a_i+ f_i|^2 d \nu \\= &\inf_{c_i\in R} \int_{\Si^2} \sum_{i=1}^3 |c_i+ x_i|^2 d \nu 
\end{align*}

For any $ (c_1, c_2, c_3) \in \R^3$, let
$$
\eta(x) =\sum_{i=1}^3 |c_i+ x_i|^2, \quad x \in \Si^2. 
$$
We know that 
$$
\eta(x) =1+ \sum_{i=1}^3 c_i^2+\sum_{i=1}^3 2c_ix_i \ge 0
$$
and
$$
\ln(\int_{\Si^2} \eta d\nu) = \ln(\int_{\Si^2} \eta f d\omega ) \ge
\int_{\Si^2} \ln{(\eta f)}  d\omega \ge \int_{\Si^2} \ln{(\eta)}  d\omega +\int_{\Si^2} \ln{f}  d\omega.
$$
We claim that  for any $ (c_1, c_2, c_3) \in \R^3$,
$$
\int_{\Si^2} \ln{(\eta)}  d\omega \ge 0.
$$
For this purpose,  we assume  without loss of generality that 
$(c_1, c_2, c_3)=(0, 0, t), t\ge 0 \in \R$ after a possible rotation,  
and define
$$
g(t)=\int_{\Si^2} \ln{(\eta)}  d\omega =\int_{\Si^2} \ln{(1+t^2+2t x_3}) d\omega=\frac{1}{2}\int_{-1}^1 \ln{(1+t^2+2t x_3} ) dx_3.
$$
Straightforward computations lead to
$$
g(t)= \ln(1+t) (1+ \frac{t^2+1}{2t})+ \ln(|t-1|)  ( \frac{t^2+1}{2t}-1 ) -1, \,\, t \in  (0, 1) \cup(1, \infty)
$$
and $g(1)=2 \ln2-1$.
It is easy to check that 
$\lim_{t \to 0^+} g(t)=0 $ and $\lim{t \to \infty} g(t)= \infty$ and
$g(t)$ is continuous in $(0, \infty)$ and $g \in C^1(0, \infty)$.
Indeed,  $g'(1)= 1 $ and
$$
g'(t)=\frac{1}{4t^2}\big( 4t+(t^2-1) \ln(\frac{1+t}{1-t})^2 \big), 
\quad t \in  (0, 1) \cup(1, \infty)
$$
is continuous in $ (0, \infty) $.  

It is easy to check by differentiation again that 
$g'(t) >0, t \in (0, \infty)$ and hence $g(t)>0,  t \in (0, \infty)$.
This proves the claim.

Hence, 
$$
\ln(D_1/D_1)  \ge \int_{\Si^2} \ln{f}  d\omega
$$
and \eqref{monotonicity} follows.
This completes the proof. 

\begin{remark}
If we choose $f=e^{2u}$, \eqref{ineq1} is equivalent to
$$
 \ln D_1 -2\int_{\Si^2} f d\omega \le \frac{4}{3} \int_{\Si^2} |\nabla u|^2 d\omega. 
 $$
 The factor $\frac{4}{3}$ in the above inequality makes the inequality weaker than the Szeg\"o limit theorem on $S^1$, where the factor is $1$ for
 the corresponding term.  Given the optimal constant in \eqref{ineq1},  we 
 may not expect that the Szeg\"o limit theorem on $\Si^2$ holds fully in its original $S^1$ form.  
 \end{remark}

$\mathbf{Acknowledgement}$  
The research was partially  done when the second author visited Princeton University,  he wishes to thank the department of mathematics 
of Princeton University for the hospitality.   Research of the first author is partially supported by NSF grant DMS- 1607091; Research of  the second author is partially supported by  NSF grants DMS-1601885 and DMS-1901914 and  Simons Foundation  Award 617072.

\bibliographystyle{plain}
\bibliography{Inequalities}
\end{document}